\documentclass[12pt]{article}
\usepackage[T2A]{fontenc}
\usepackage[utf8]{inputenc}
\usepackage[english]{babel}
\usepackage{amsmath,amsfonts,amssymb}
\usepackage[mathscr]{eucal}
\usepackage{graphicx}  
\usepackage{icomma}
\usepackage{bbm}
\usepackage{tikz}
\usetikzlibrary{shapes,snakes}
\usepackage{color}
\usepackage{algorithmic}
\usepackage{verbatim}
\usepackage{xurl}
\usepackage{subfig}
\usepackage{array}

\DeclareMathOperator{\rank}{rank}

\def\Re{{\rm Re\,}}

\makeatletter
\renewcommand*{\@biblabel}[1]{\hfill#1.}
\makeatother

\usepackage[T2A]{fontenc}
\usepackage{amsmath,amsfonts,amssymb}
\usepackage[mathscr]{eucal}
\usepackage{graphicx}  
\usepackage{icomma}
\usepackage{bbm}
\usepackage{tikz}
\usetikzlibrary{shapes,snakes}
\usepackage{color}
\usepackage{algorithmic}
\usepackage{verbatim}
\usepackage{xurl}
\usepackage{subfig}
\usepackage{array}

\usepackage{amsthm}

\def\Re{{\rm Re\,}}

\DeclareMathOperator*{\argmax}{arg\,max}

\DeclareMathOperator*{\sign}{sign}

\makeatletter
\renewcommand*{\@biblabel}[1]{\hfill#1.}
\makeatother

\newtheorem{theorem}{Theorem}
\newtheorem{lemma}{Lemma}
\newtheorem{statement}{Statement}

\begin{document}
\vspace*{0mm}

% УДК  XXX.XX
\begin{center} 
\large\bf On the algorithm of best approximation \\
\vspace*{1mm}
by low rank matrices in the Chebyshev norm.\footnote{$^)$\,This work was supported by Russian Science Foundation (project 21-71-10072).}$^)$
\end{center} 

\vspace*{1mm}
\begin{center}
{\large\rm\bf
% \copyright\,\,2021 г.\,\,\,\,\,\,\,  
S. Morozov$^{1*}$, N. Zamarashkin$^{1**}$, E. Tyrtyshnikov$^{1***}$}\\
\vspace*{3mm}
{\it 
	$^{1}$ 
	Marchuk Institute of Numerical Mathematics of the Russian Academy of Sciences\\
e--mail: $^{*}$stanis-morozov@yandex.ru\\
e--mail: $^{**}$nikolai.zamarashkin@gmail.com\\
e--mail: $^{***}$eugene.tyrtyshnikov@gmail.com
} \\
\medskip
% Received\,\, \,\,\,\,--.--.2021

% Revised\,\, \,\,\,--.--.2021

% Accepted\,\, \,\,\,\,--.--.2021
\vspace*{2mm}
\end{center}

\noindent
{\small The low-rank matrix approximation problem is ubiquitous in computational mathematics. Traditionally, this problem is solved in spectral or Frobenius norms, where the accuracy of the approximation is related to the rate of decrease of the singular values of the matrix. However, recent results indicate that this requirement is not necessary for other norms. In this paper, we propose a method for solving the low-rank approximation problem in the Chebyshev norm, which is capable of efficiently constructing accurate approximations for matrices, whose singular values do not decrease or decrease slowly.% Библ. 12. Фиг. 3.
}

\vspace*{2mm}


\vspace*{3mm}

\begin{center}
{\large 1. INTRODUCTION}
\end{center}
Low-rank matrices are ubiquitous in science. They serve as a tool for low-parametric matrix approximation in numerous applications such as computational mathematics \cite{bebendorf2008means}, computational fluid dynamics \cite{son2014data}, recommender systems \cite{he2016fast}, machine learning \cite{yang2018oboe}, and others.

However, one typically assumes that the singular values of the matrix that needs to be approximated decay rapidly. This assumption is made, primarily, because there are efficient algorithms for close-to-optimal low-rank approximation in unitarily invariant norms \cite{GTZ-1997, HMT-2011, OZ-2018}.

On the other hand, in modern applications, especially in the field of big data, it is often more natural to use other matrix norms. For example, in the classical scheme of recommender systems, one deals with rating matrices, the rows of which correspond to products (films, music, etc.) and the columns correspond to users. The values of the matrix entries determine the ratings given to products by users. To restore the missing ratings and provide recommendations, a low-rank approximation of the matrix is constructed based on the known elements. In this case, it seems more natural to approximate the matrix elementwise, rather than in spectral or Frobenius norms: try to approximate the values of all ratings in the best possible way. Moreover, as follows from the article \cite{udell2019big}, statistical models used to describe the rating matrices lead, in general, to matrices with a slow decay of singular values, but which can be approximated elementwise with matrices of low rank. The latter means that the use of low-rank approximation algorithms based on singular value decompositions in recommender systems can hardly be considered reasonable.

We consider the best low-rank appromation problem in the so-called Chebyshev norm:
\begin{equation*}
% \label{cheb_norm}
\|X\|_C = \max\limits_{i,j} |x_{ij}|.
\end{equation*}
Namely, for given matrix $ A \in \mathbb{C}^{m \times n}$ and integer $r$, it is required to find $ U \in \mathbb {C}^{m \times r} $ and
$ V \in \mathbb {C}^{n \times r} $ such that
\begin{equation}
\label{global_problem}
\mu = \inf\limits_{U \in \mathbb{C}^{m \times r}, V \in \mathbb{C}^{n \times r}} \left\| A - UV^T \right\|_{C}.
\end{equation}
Matrices $ \widehat{U} $ and $ \widehat{V} $ satisfying 
\begin{equation*}
    \left\| A - \widehat{U}\widehat{V}^T \right\|_{C} = \mu,
\end{equation*}
will be called matrices of best approximation of $A$ of rank $r$.

As natural as problem (\ref{global_problem}) is, to date, it has been little studied: there are asymptotic estimates for the accuracy of the approximation (\ref{global_problem}) \cite{udell2019big} and a method for constructing local minima of the problem (\ref{global_problem}) in the case of rank 1 \cite{daugavet1971}.

In this paper, we propose and justify an algorithm for solving the problem
\begin{equation*}
% \label{one_matrix_problem}
\mu = \inf\limits_{U \in \mathbb{R}^{m \times r}} \left\| A - UV^T \right\|_{C}
\end{equation*}
for an arbitrary rank. Based on this algorithm, we develop a method of finding local minima of the problem (\ref{global_problem}) for an arbitrary rank. A large number of numerical experiments show that the asymptotic estimates proved in \cite{udell2019big} are generally not optimal.

The rest of the paper is organized as follows. In Section~2, we collect the known results from the literature about the problem we consider. Section~3 presents basic results on the properties of the solution, including the questions of the existence, uniqueness, and continuity of the solution. In addition, we discuss the question of existence of characteristic sets and their properties, together with the well-known results on how to solve the problem for $ (r + 1) \times r $ matrices. We also consider the optimality criteria of solutions. In Section~4, we present a combinatorial formula for solving the problem and propose a generalized Remez algorithm, which makes it possible to find solutions, in practice, in a polynomial number of operations. Section~5 gives an algorithm for solving the problem in the case when both matrices $ U $ and $ V $ are assumed to be unknown. Numerical experiments from Section~6 demonstrate the effectiveness of the proposed method, and also lead to a number of new questions about the asymptotic accuracy of approximations of matrices in the Chebyshev norm. Section~7 concludes the paper.

\begin{center}
{\large 2. EXISTING RESULTS}
\end{center}
As far as we know, the problem of constructing and analyzing low-rank approximations of matrices in the Chebyshev norm has been little studied. We rely on two works \cite{udell2019big, daugavet1971}. The first of them contains results on the asymptotic properties of Chebyshev approximations of matrices (without the assumption that singular values decrease), and the second contains a method for finding local optima for the problem (\ref{global_problem}) for rank $1$.

Let us take a closer look at these works. One of the results proven in \cite{udell2019big} is the following general statement.
\begin{theorem}
\label{townsend_theorem}
Let $X \in \mathbb{R}^{m\times n}$, where $m \ge n$ and $0 < \varepsilon < 1$. Then, with
\begin{equation*}
r = \lceil 72\log{(2n + 1)}/\varepsilon^2 \rceil
\end{equation*}
we have
\begin{equation*}
\inf\limits_{\rank Y \le r} \|X - Y\|_C \le \varepsilon \|X\|_2.
\end{equation*}
\end{theorem}
It is clear from this theorem that Chebyshev low-rank approximations have great potential. So for any sequence of matrices with bounded spectral norms, for a fixed accuracy $\varepsilon$, the rank of the Chebyshev approximation grows at most logarithmically. For example, for the best rank$-(n-1)$ approximation of the $n \times n$ identity matrix in the spectral or Frobenius norms, its accuracy is $ 1 $. At the same time, in the Chebyshev norm, the identity matrix can be approximated with any fixed accuracy $ \varepsilon > 0 $ with rank that is logarithmic in the order of the matrix $ n $. In \cite{udell2019big}, this property of the Chebyshev norm is called one of the main reasons why matrices arising in data analysis can be efficiently approximated by low-rank matrices.

As far as we know, the only paper that studies the problem (\ref{global_problem}) is \cite{daugavet1971}, where rank-$1$ approximation is considered. First, in \cite{daugavet1971}, the problem of the form 
\begin{equation}
\label{one_vector_problem}
\mu = \inf\limits_{u \in \mathbb{R}^{m}} \left\| A - uv^T \right\|_{C}.
\end{equation}
is solved. It is easy to see that for each row of the matrix $ A $ the problem (\ref{one_vector_problem}) can be solved independently and, therefore, is reduced to
\begin{equation*}
% \label{one_number_problem}
\mu = \inf\limits_{u \in \mathbb{R}} \left\| a - uv \right\|_{\infty}.
\end{equation*}
This implies a simple algorithm for solving the problem (\ref{one_vector_problem}). In order to obtain the local minimum of the solution to the problem (\ref{global_problem}), the authors use the alternance method. Let an initial vector $v^{(0)}$ be given. Solving the problem (\ref{one_vector_problem}) for fixed $v=v^{(0)}$, we find a solution $u^{(1)}$. Further, for fixed $ u = u^{(1)} $, we solve the problem
\begin{equation*}
\mu = \inf\limits_{v \in \mathbb{R}^{n}} \left\| A - uv^T \right\|_{C}
\end{equation*}
and find a solution $v^{(1)}$. Continuing according to this scheme, we come to a solution, which, however, is not always a global optimum to the problem (\ref{global_problem}). In addition, \cite{daugavet1971} contains a necessary condition for the optimality of the solution $ u, v $ of the problem (\ref{global_problem}) for rank 1. For simplicity, we assume that all elements $ u $ and $ v $ are nonzero.

\begin{statement}
\label{simple_necessary_condition}
Let all elements of the vectors $ u $ and $ v $ be nonzero and assume that they are a solution to the problem (\ref{global_problem}). Let $ R = A - uv^T $. Then in the matrix $ R $ there is a cycle, that is, a set of indices $(i_1, j_1), (i_1, j_2), (i_2, j_2), \dots, (i_k, j_k), (i_k, j_1)$ such that
\begin{enumerate}
\item indices $i_1, \dots, i_k$ are distinct;
\item indices $j_1, \dots, j_k$ are distinct;
\item at each of these positions in the matrix $ R $, the maximum absolute value is reached;
\item let $(i_t, j_p)$ and $(i_g, j_h)$  be adjacent in the cycle, that is, they are different pairs of indices such that $t=g$ or $p=h$. Then the signs of the values $u_{i_t} v_{j_p} r_{i_t j_p}$ and $u_{i_g} v_{j_h} r_{i_g j_h}$ are distinct.
\end{enumerate}
\end{statement}
In what follows, we prove generalizations of these results to the case of arbitrary rank.

% \section{Предварительные результаты}
\begin{center}
{\large 3. PRELIMINARIES}
\end{center}
Now, we are going to provide several basic results that will be useful to us. A significant part of them is a reformulation of the well-known results from the theory of Chebyshev approximations of functions for the matrix case \cite{dzyadyk1977book, smirnov1964constructive}. In this section, we are interested in the following problem
\begin{equation}
\label{one_matrix_problem_statement}
    \mu = \inf\limits_{U \in \mathbb{C}^{m \times r}} \left\| A - UV^T \right\|_{C}.
\end{equation}
Matrix $\widehat{U}$ satisfying
\begin{equation*}
    \left\| A - \widehat{U}V^T \right\|_{C} = \mu,
\end{equation*}
is called a matrix of the best approximation $ A $ with respect to the system of vectors $ V $. The task is to find a matrix $ \widehat{U} \in \mathbb{C}^{m \times r} $ of the best approximation. It is easy to understand that this problem is divided into $ m $ independent subproblems for each row of the matrix $ A $, for which it is required to find the corresponding row of the matrix $ \widehat{U} $. Therefore, we are going to solve the following problem. Let the matrix $ V \in \mathbb{C}^{n \times r} $ and the vector $ a \in \mathbb{C}^{n} $ be given. It is required to find
\begin{equation}
\label{vector_problem_statement}
    \mu = \inf\limits_{u \in \mathbb{C}^{r}} \left\| a - Vu \right\|_{\infty}
\end{equation}
and vector $\widehat{u} \in \mathbb{C}^r$ satisfying
\begin{equation*}
    \left\| a - V\widehat{u} \right\|_{\infty} = \mu.
\end{equation*}

% \subsection{Существование, единственность, непрерывность}
\begin{center}
	{\it 3.1. Existence, uniqueness, continuity}
\end{center}

% Изучим вопрос 
Let us present the results on the existence, uniqueness and continuity of the solution to the problem (\ref{vector_problem_statement}). The existence of a solution, that is, the existence of such $ \widehat{u} \in \mathbb{C}^r $ that
\begin{equation*}
    \left\| a - V\widehat{u} \right\|_{\infty} = \mu
\end{equation*}
is obvious. Consider the question of the uniqueness of the solution. In what follows, we will denote columns of matrices by subscripts and rows by superscripts. Let us introduce the concept of a Chebyshev vector system.

{\bf Definition 1.} The columns of the matrix $ V \in \mathbb{C}^{n \times r} $ form a Chebyshev system of vectors if any $ r $ rows of the matrix $ V $ are linearly independent. 

This concept is closely related to the uniqueness of the solution to the problem of best approximation. Namely, we have the following 
\begin{theorem}
(Haar \cite{dzyadyk1977book}). Let the matrix $ V \in \mathbb{C}^{n \times r} $ and $ n> r $ be given. Then, for any vector $ a \in \mathbb{C}^n $ the solution is unique if and only if the columns of the matrix $ V $ form a Chebyshev system.
\end{theorem}
The following statement is also true
\begin{statement}
\label{optimal_solution_has_equal_elements}
Let the columns of the matrix $ V \in \mathbb {C}^{n \times r} $, where $ n> r $ form a Chebyshev system and let $ V \widehat{u} $ be the vector of the best approximation. Then the maximum absolute value is reached at least at $ r + 1 $ points, that is, there are $ i_1, \dots, i_{r + 1} $ which satisfy the equality
\begin{equation*}
    |a_{i_j} - (V\widehat{u})_{i_j}| = \|a - V\widehat{u}\|_\infty, ~~~ j = 1, \dots, r + 1.
\end{equation*}
\begin{proof}
Let the number of points at which the maximum absolute value is reached be equal to $ r_1 <r + 1 $. Then, having solved the system with $r_1$ equations and $ r $ unknowns, whose rows are linearly independent, we can obtain a vector $ p \in \mathbb {C} ^ r $ such that
\begin{equation*}
    (Vp)_{i_j} = a_{i_j} - (V\widehat{u})_{i_j}, ~~~ j = 1, \dots, r_1.
\end{equation*}
But then the vector $ V (\widehat{u} + \delta p) $ for a sufficiently small $ \delta $ deviates from $ a $ less than $ V \widehat{u} $. We come to a contradiction.
\end{proof}
\end{statement}
Consider the issue of continuity of the solution.
\begin{theorem}
\label{Nikolskiy_continuity}
(Nikolsky \cite{dzyadyk1977book}). Let the system of columns of the matrix $ V \in \mathbb {C}^{n \times r} $, where $ n> r $, form a Chebyshev system. Then the coefficients of the best approximation vector $ \widehat{u} $ continuously depend on the approximated vector $ a $, and the column system $ V $, that is, $ \forall \varepsilon> 0 $ $ \exists \delta = \delta (a, V , \varepsilon)> 0 $ such that if $ \| a - b \|_\infty + \| V - W \| <\delta $, then $ \| \widehat{u}(a, V) - \widehat{u} (b, W) \|_\infty < \varepsilon $, where $ \widehat{u} (a , V) $ and $ \widehat{u}(b, W) $ denote the coefficients of the optimal solution for the vectors $ a $ and $b$ by the systems $ V $ and $ W $, respectively.
\end{theorem}

% \subsection{Характеристические множества}
\begin{center}
	{\it 3.2. Characteristic sets }
\end{center}

Let $ J $ denote the set of indices $ J = \left \{1,2, \dots, n \right \}, $ and let $ J'$ and $ J''$ be subsets of $ J $. Denote
\begin{equation*}
\mu(J') = \inf\limits_{u \in \mathbb{C}^{r}} \left\| a(J') - V(J') u \right\|_{\infty},
\end{equation*}
where $ V(J') $ denotes the submatrix of the matrix $ V $ containing rows with numbers from the set $ J' $, and $ a(J') $ denotes the subvector of the vector $ a $ containing elements with the indices from $ J' $.

{\bf Definition 2.} A set $ J' $ is called a characteristic set if $ \mu(J) = \mu(J') $ and for any subset  $J'' \subsetneq J'$, $\mu(J'') < \mu(J)$.

Further, we will show that if the columns of $ V $ are linearly independent, then there is at least one characteristic set containing at most $ 2r + 1 $ points in the complex case and at most $ r + 1 $ points in the real case.

We need the following notation \cite{dzyadyk1977book, smirnov1964constructive}. Let $\lambda \ge 0$.
Denote
\begin{eqnarray*}
	F(j, u) & =  & \left| a_{j} - u^T v^j \right|, \\
        K(j, \lambda) & =  & \left\{ u \in \mathbb{C}^{r} | F(j, u) \le \lambda \right\}, \\
        K(J', \lambda) & = & \bigcap\limits_{j \in J'} K(j, \lambda) = \left\{ u \in \mathbb{C}^{r} | F(j, u) \le \lambda, \forall j \in J' \right\}.
\end{eqnarray*}
It is easy to see that the following is true.
\begin{statement}
\begin{eqnarray*}
	K(j, \lambda') & \subset & K(j, \lambda''), ~~~ K(J', \lambda') \subset K(J', \lambda''), ~~~ 0 \le \lambda' < \lambda'', \\
	K(J'', \lambda) & \subset &  K(J', \lambda), ~~~ J' \subset J''.
\end{eqnarray*}
\end{statement}

\begin{lemma}
\label{existance_lemma_matrix}
Let the columns of the matrix $ V $ be linearly independent. Then the set $ K (j, \lambda) $ is convex and closed, and the set $ K (J, \lambda) $ is bounded for any $ \lambda \ge 0 $.
\begin{proof}
Let us prove that the sets $ K (j, \lambda) $ are closed. Let $ u_1 $ be the limit point of the set $ K (j, \lambda) $. Then in any of its neighborhood there are points $ u \in K (j, \lambda) $.
\begin{equation*}
    F(j, u_1) \le |F(j, u_1) - F(j, u)| + F(j, u)
\end{equation*}
Since $u \in K(j, \lambda)$, then $F(j, u) \le \lambda$. The function
\begin{equation*}
    F(j, u) = \left| a_j - u^T v^j \right|
\end{equation*}
is continuous in $ u $ for any fixed $ j $. Then for any $\varepsilon > 0$ there is $\delta > 0$ such that if $|u - u_1| < \delta$, then $|F(j, u_1) - F(j, u)| < \varepsilon$. Thus, we have
\begin{equation*}
    F(j, u_1) < \varepsilon + \lambda
\end{equation*}
for any $\varepsilon > 0$, hence $u_1 \in K(j, \lambda)$. The closedness is proved.

Let us prove that $ K (J, \lambda) $ is bounded. Consider the vector $ Vu $ for $ \| u \|_1 = 1 $. The quantity $ \| Vu \|_\infty $ is a function continuous in $ u $ on a compact set; therefore, at some point it reaches its minimum value 
\begin{equation*}
    M = \|V\widehat{u}\|_\infty \le \|Vu\|_\infty.
\end{equation*}
Since the columns $ V $ are linearly independent, any non-trivial linear combination of them is not equal to $ 0 $ and $ M> 0 $. 

Let $\|u\|_1 \ge \dfrac{C + 1}{M}$, where $C > 0$ is a constant.
Then
\begin{equation*}
    \|a - Vu\|_\infty \ge \|Vu\|_\infty - \|a\|_\infty \ge \|u\|_1 M - \|a\|_\infty
    \ge C + 1 - \|a\|_\infty.
\end{equation*}
Then for $\|u\|_1 \ge \dfrac{C + 1}{M}$ the condition $F(J, u) \le \lambda = C - \|a\|_\infty$ cannot be satisfied, that is $u \notin K(j, \lambda)$ for $\lambda \le C - \|a\|_\infty$. Since $ C $ is arbitrary, the boundedness is proved.

Let us prove the convexity of the set $K(j, \lambda)$. Let $u_1, u_2 \in K(j, \lambda)$, that is $F(j, u_1) \le \lambda$ and $F(j, u_2) \le \lambda$. We need to prove that $F(j, \tau u_1 + (1 - \tau) u_2) \le \lambda$ for any $\tau \in (0, 1)$.
\begin{eqnarray*}
	F(j, \tau u_1 + (1 - \tau)u_2) & = & | a_j - (\tau u_1 + (1 - \tau)u_2)^Tv^j | \\ 
	& = & 
    | \tau(a_j - u_1^T v^j) + (1-\tau)(a_j - u_2^T v^j) | \\ 
       & \le &  \tau F(j, u_1) + (1 - \tau)F(j, u_2) \le \lambda.
\end{eqnarray*}
\end{proof}
\end{lemma}

We denote by $ M_k $ the set of all ordered subsets $ J_k $, consisting of $ k $ elements $ i_1, \dots, i_k $, taken from the set $ J $. We denote by $ \mu_k (J) $ the exact upper bound of the least deviations from zero of the function $ F (j, u) $ on all subsets of $ J_k \in M_k $: 
\begin{equation*}
    \mu_k(J) = \max\limits_{J_k \in M_k} \mu(J_k) = \max\limits_{J_k \in M_k} \min\limits_{u \in \mathbb{C}^r} \max\limits_{j \in J_k} F(j, u).
\end{equation*}
It is easy to prove the following
\begin{statement}
$\mu_k(J) \le \mu_{k+1}(J) \le \mu(J), ~~~ \forall k$.
\end{statement}

For further analysis, we need the following theorem 
\begin{theorem}
(Helly). If the collection $ K $ of closed convex sets of points $ x \in \mathbb{R}^r $ contains at least $ r + 1 $ sets (some of which may be the same), the intersection of any $ r + 1 $ sets from $ K $ is not empty and the intersection of a finite number of sets from $ K $ is bounded, then the intersection of all sets from $ K $ is not empty. 
\end{theorem}

Based on Helly's theorem, let us prove the following result.
\begin{theorem}
\label{schnirelman_matrix}
(Shnirelman \cite{dzyadyk1977book}). If there is a $\lambda_0 > \mu_{2r+1}$ ($\lambda_0 > \mu_{r+1}$ in real case), such that for any $j$ and any $\lambda$, $\mu_{2r+1} < \lambda < \lambda_0$ ($\mu_{r+1} < \lambda < \lambda_0$ in real case), the sets $K(j, \lambda)$ are closed and convex and the intersection of a finite number of sets $K(j, \lambda)$ is bounded, then $\mu_{2r+1}(J) = \mu(J)$ ($\mu_{r+1}(J) = \mu(J)$ in real case).
\begin{proof}
Let $ k = r + 1 $ in the real case and $ k = 2r + 1 $ in the complex case. Since the sets $ J_k $ and $ M_k $ are finite, we have that the maxima and minima are attained, therefore, for any $ \lambda > \mu_k (J) $, the set $ K(J_k, \lambda) $ is not empty. Moreover, since $ K(j, \lambda) $ are convex and closed, then 
\begin{equation*}
    K(J_{k}, \lambda) = \bigcap\limits_{j \in J_{k}} K(j, \lambda)
\end{equation*}
is nonempty and convex for any $J_{k} \in M_{k}$. 

Let us verify that all conditions of Helly's theorem are satisfied. As the collection of sets $ K $ we take the sets $ K (j, \ lambda) $, $ j \ in J $. By the hypothesis of the theorem, they are closed and convex, and the intersection of some finite number of these sets is bounded. In the real case, the fact that the intersection of any $ j + 1 $ sets is not empty is equivalent to the fact that $ K(J_{r + 1}, \lambda) $ is not empty, as shown above. In the complex case, we need to work with the space $ \mathbb{C}^r $, which we identify with $ \mathbb{R}^{2r} $, so we need the intersection of any $ 2r + 1 $ sets to be non-empty , which was also shown above. So, all conditions of Helly's theorem are satisfied and we have that 
\begin{equation*}
    K(J, \lambda) = \bigcap\limits_{j \in J} K(j, \lambda)
\end{equation*}
is nonempty, convex and bounded.

Assume that the sequence $\{\lambda_t\}$ decreases, $\mu_k < \lambda_t < \lambda_0$ and it tends to $\mu_k$. Moreover, let $K(J, \lambda_{t+1}) \subset K(J, \lambda_t)$ and the intersection
\begin{equation*}
    K = \bigcap\limits_{t=1}^{\infty} K(J, \lambda_t)
\end{equation*}
be nonempty, convex and bounded. Let $u_0 \in K$.
This means that $ u_0 \in K(J, \lambda_t) $, hence $ F (j, u_0) \le \lambda_t $ for any $ j \in J $ and any $ t $. In the limit as $ t \to \infty $, we see that $ F (j, u_0) \ le \mu_k(J) $ for any $ j \in J $, whence
\begin{equation*}
    \mu(J) \le \max\limits_{j \in J} F(j, u_0) \le \mu_k(J).
\end{equation*}
But, as noted above, $\mu_{k}(J) \le \mu(J)$.
\end{proof}
\end{theorem}

This theorem allows us to formulate the following result.
\begin{theorem}
Let the columns of the matrix $ V \in \mathbb{C}^{n \times r} $, where $ n \ge r $, be linearly independent and assume that the vector $ a $ does not belong to the range of the matrix $ V $. Then there is at least one characteristic set consisting of at most $ 2r + 1 $ points in the complex case and $ r + 1 $ points in the real one. Moreover, if the system of columns of the matrix $ V $ is Chebyshev, then any characteristic set consists of at least $ r + 1 $ points. 
\begin{proof}
The result for an arbitrary system immediately follows from Lemma~\ref{existance_lemma_matrix} and Theorem~\ref{schnirelman_matrix}. In the case of a Chebyshev system, on any set of $ r $ or fewer points, one can solve the system and precisely approximate the vector at these points, and since the set is characteristic, this contradicts the condition that $ a $ does not belong to the range of $ V $. 
\end{proof}
\end{theorem}

%\begin{note}
%В комплексном случае все значения от $r+1$ до $2r + 1$ достижимы. См. Раздел~\ref{complex_case_notes}.
%\end{note}

% \subsection{О задаче поиска равноудаленных точек}
\begin{center}
	{\it 3.3. On the problem of finding equidistant points }
\end{center}

Let us introduce the concept of an equidistant point of the system

{\bf Definition 3.}
Let $V \in \mathbb{R}^{(r+1)\times r}$ and $a \in \mathbb{R}^{r + 1}$. Let the system
\begin{equation*}
    Vu = a
\end{equation*}
be inconsistent. A point $ u $ is called an equidistant point of the system if 
\begin{equation*}
    \rho(u) = |(v^1, u) - a_1| = |(v^2, u) - a_2| = \dots = |(v^{r+1}, u) - a_{r+1}|.
\end{equation*}
A point $ u $ is called the best equidistant point of the system if it is equidistant and the value $ \rho(u) $ is minimal. 

Let us present the results \cite{dzyadyk1974uniform} on the structure of the set of all equidistant points of the system in the real case. Let
\begin{equation*}
    \widehat{V}_j = \begin{bmatrix}
        v_1^1 & v_2^1 & \dots & v_r^1 \\
        v_1^2 & v_2^2 & \dots & v_r^2 \\
        \vdots & \vdots & \vdots & \vdots \\
        v_1^{j-1} & v_2^{j-1} & \dots & v_r^{j-1} \\
        v_1^{j+1} & v_2^{j+1} & \dots & v_r^{j+1} \\
        \vdots & \vdots & \vdots & \vdots \\
        v_1^{r+1} & v_2^{r+1} & \dots & v_r^{r+1} \\
    \end{bmatrix}
\end{equation*}
and $\widehat{a}_j = (a_1, a_2, \dots, a_{j-1}, a_{j+1}, \dots, a_{r+1})^T$. Then denote $D_j = \det \widehat{V}_j$ and let $\widehat{u}^j$ be the solution of the system $\widehat{V}_j u = \widehat{a}_j$. The following theorems about the set of equidistant points of the system hold.
\begin{theorem}
\label{full_equidistant_solution}
(Dzyadyk \cite{dzyadyk1974uniform}). Let an inconsistent system of equations $Vu = a$ be given, where $V \in \mathbb{R}^{(r+1)\times r}$ and $a \in \mathbb{R}^{r + 1}$. Then
\begin{enumerate}
    \item For each $j=1,\dots,r+1$ we have
    \begin{equation*}
        (v^j, \widehat{u}^j) - a_j = \dfrac{(-1)^{j+1}}{D_j} \sum\limits_{\nu=1}^{r+1} (-1)^{\nu} a_{\nu} D_{\nu}.
    \end{equation*}
    \item For any real $k_j$, $j=1,\dots,r+1$ such that
    \begin{equation*}
        \sum\limits_{j=1}^{n+1} |D_j| e^{i k_j} \neq 0
    \end{equation*}
    point $ u $ determined by the formula 
    \begin{equation*}
        u = \rho \sum\limits_{j=1}^{r+1} \dfrac{\widehat{u}^j e^{i k_j}}{|(v^j, \widehat{u}^j) - a_j|} = \rho \dfrac{\sum\limits_{j=1}^{r+1} |D_j| \widehat{u}^j e^{i k_j} }{\left| \sum\limits_{j=1}^{r+1} (-1)^j D_j a_j \right|},
    \end{equation*}
    where
    \begin{equation*}
        \rho = \left( \sum\limits_{j=1}^{r+1} \dfrac{e^{i k_j}}{|(v^j, \widehat{u}^j) - a_j|} \right)^{-1} = \dfrac{\left| \sum\limits_{j=1}^{r+1} (-1)^j D_j a_j \right|}{ \sum\limits_{j=1}^{r+1} |D_j| e^{i k_j}}
    \end{equation*}
    is an equidistant point of the system $ Vu = a $, while $ |\rho| $ is equal to the $V$-distance from the point $ u $ to $ a $. 
    \item Any equidistant point $ u $ of the system $ Vu = a $ can be represented for some real $ k_j $ by the formula above. In this case, $ k_j $ can, in particular, be expressed by the formula 
    \begin{equation*}
        k_j = \arg ((v^j, u) - a_j) - \arg((v^j, \widehat{u}^j) - a_j).
    \end{equation*}
\end{enumerate}

\end{theorem}

\begin{theorem}
\label{best_equidistant_solution}
(Dzyadyk \cite{dzyadyk1974uniform}). Let an inconsistent system of equations $Vu = a$ be given, where $V \in \mathbb{R}^{(r+1)\times r}$ and $a \in \mathbb{R}^{r + 1}$. Then the best equidistant point $u^*$ of the system $Vu = a$ can be determined by the formula 
    \begin{equation*}
        u^* = \rho^* \sum\limits_{j=1}^{r+1} \dfrac{\widehat{u}^j}{|(v^j, \widehat{u}^j) - a_j|} = \dfrac{\sum\limits_{j=1}^{r+1} |D_j| \widehat{u}^j}{\sum\limits_{j=1}^{r+1} |D_j|},
    \end{equation*}
    where
    \begin{equation*}
        \rho^* = \left( \sum\limits_{j=1}^{r+1} \dfrac{1}{|(v^j, \widehat{u}^j) - a_j|} \right)^{-1} 
	    = \dfrac{\left| \sum\limits_{j=1}^{r+1} (-1)^j D_j a_j \right|}{ \sum\limits_{j=1}^{r+1} |D_j|}.
    \end{equation*}
\begin{proof}
It suffices to note that the value
\begin{equation*}
    \rho = \left( \sum\limits_{j=1}^{r+1} \dfrac{e^{i k_j}}{|(v^j, \widehat{u}^j) - a_j|} \right)^{-1}
\end{equation*}
takes the smallest value when all terms are coaligned, that is, $e^{i k_1} = e^{i k_2} = \dots = e^{i k_{r+1}}$.
\end{proof}
\end{theorem}

% \subsection{Критерии оптимальности}
\begin{center}
	{\it 3.4. Optimality criteria}
\end{center}

Let us give several criteria for the optimality of the solution. These criteria are interesting themselves, and allow one to get new important information about the problem.
\begin{theorem}
(Kolmogorov). 
Let a system of vectors with a matrix $ V \in \mathbb{C}^{n \times r} $ and a vector $ a \in \mathbb{C}^n $ be given, so that vector $a$ should be approximated by a linear combination of columns of $ V $. For the vector $ V \widehat {u} $ to be the vector of the best approximation for $ a $, it is necessary and sufficient that on the set $ E = E (V \widehat {u}) $ of all points at which for the vector $ V \widehat{u} $ the maximum absolute value of the residual is reached, for all vectors of the form $ Vu $ the equality
\begin{equation*}
    \min\limits_{j \in E} \Re \left( (Vu)_j \overline{ \left( a_j - (V\widehat{u})_j \right)} \right) \le 0.
\end{equation*}
is satisfied.
\begin{proof}
Necessity. Let $ V \widehat{u} $ be the vector of best approximation for $ a $. By contradiction, let 
\begin{equation*}
    \min\limits_{j \in E} \Re \left( (Vu)_j \overline{ \left( a_j - (V\widehat{u})_j \right)} \right) > c > 0
\end{equation*}
for some vector $u \in \mathbb{C}^r$. Denote
\begin{equation*}
    G = \max\limits_{j \in E} |a_j - (V\widehat{u})_j|, ~~~ G' = \max\limits_{j \notin E} |a_j - (V\widehat{u})_j|,
\end{equation*}
\begin{equation*}
    h = G - G' > 0, ~~~ M = \max\limits_{j} |(Vu)_j|, ~~~ \lambda = \max\left\{ \dfrac{c}{M^2}, \dfrac{h}{2M} \right\} > 0.
\end{equation*}
Let us prove then that the vector $ V (\widehat{u} + \lambda u) $ approximates the vector $ a $ better. 

1) Let $j \in E$. Then
\begin{eqnarray*}
	|a_j - (\widehat{u} + \lambda u)^T v^j|^2 & = & (a_j - \widehat{u}^T v^j - \lambda u^T v^j) \cdot 
	( \overline{(a_j - \widehat{u}^T v^j)} - \lambda \overline{(u^T v^j)}) \\
	& = & |a_j - \widehat{u}^T v^j|^2 + \lambda^2 |u^T v^j|^2 - 
        2 \lambda \Re \left( (Vu)_j \overline{ \left( a_j - (V\widehat{u})_j \right)} \right) \\ 
	& \le &  
    G^2 + \lambda^2 M^2 - 2 \lambda \Re \left( (Vu)_j \overline{ \left( a_j - (V\widehat{u})_j \right)} \right) \\ 
	& < &
    G^2 + \lambda^2 M^2 - 2\lambda c \le
    G^2 + \lambda \dfrac{c}{M^2} M^2 - 2\lambda c \\ 
    & = & G^2 - \lambda c < G^2.
\end{eqnarray*}
2) Let $j \notin E$. Then
\begin{equation*}
    |a_j - (\widehat{u} + \lambda u)^T v^j| \le |a_j - \widehat{u}^T v^j| + \lambda |u^T v^j| \le G' + \lambda M \le G - h + \dfrac{h}{2M} M = G - h / 2 < G.
\end{equation*}
It follows that the vector $ V (\widehat{u} + \lambda u) $ approximates the vector $ a $ better. We get a contradiction. 

Sufficiency. Let the condition of the Kolmogorov criterion be satisfied with the vector of coefficients $ \widehat{u} $ and let $ u \in \mathbb{C}^r $ be arbitrary. 
Consider the vector $ w = V(u - \widehat{u}) $. Let us choose the index $ j_0 $ for which the inequality
\begin{equation*}
    \Re \left( (V(u - \widehat{u}))_{j_0} \overline{ \left( a_{j_0} - (V\widehat{u})_{j_0} \right)} \right) \le 0
\end{equation*}
holds. Then
\begin{eqnarray*}
	|a_{j_0} - (Vu)_{j_0}|^2 & = & |a_{j_0} - (V\widehat{u})_{j_0} - ((Vu)_{j_0} - (V\widehat{u})_{j_0})|^2 \\  
	& =  &  |a_{j_0} - (V\widehat{u})_{j_0}|^2 + |(Vu)_{j_0} - (V\widehat{u})_{j_0}|^2
    - 2 \Re \left( (V(u - \widehat{u}))_{j_0} \overline{ \left( a_{j_0} - (V\widehat{u})_{j_0} \right)} \right) \\
	& \ge &  |a_{j_0} - (V\widehat{u})_{j_0}|^2.
\end{eqnarray*}
This shows that for any vector $ u \in \mathbb {C}^r $, the approximation given by the vector $ \widehat{u} $ is not worse; that is, it is optimal.
\end{proof}
\end{theorem}

Let us give another, in some situations more convenient, optimality criterion. In a sense, it is a reformulation of the Kolmogorov criterion using the following lemma. 

\begin{lemma}
Let $u_{ij}$, $i=1, \dots, m$, $j=1, \dots, n$ be some numbers. Then the numbers $\delta_i \ge 0$, $i=1, \dots, m$ that are not all zero and such that
\begin{equation*}
    \sum\limits_{i=1}^{m} \delta_i u_{ij} = 0, ~~~ j=1, \dots, n,
\end{equation*}
exist if and only if for any system of numbers $c_j$, $j=1, \dots, n$, inequalities 
\begin{equation*}
    \Re \sum\limits_{j=1}^{n} c_j u_{ij} > 0, ~~~ i=1,\dots,m,
\end{equation*}
are not satisfied simultaneously.
\begin{proof}
Necessity. Assume that for some $\delta_i \ge 0$, $i=1,\dots,m$
\begin{equation*}
    \sum\limits_{i=1}^{m} \delta_i u_{ij} = 0, ~~~ j=1,\dots,n.
\end{equation*}
Then
\begin{equation*}
    \sum\limits_{i=1}^{m} \delta_i \Re \sum\limits_{j=1}^n c_j u_{ij} = \Re \sum\limits_{j=1}^{n}c_j \sum\limits_{i=1}^{m} \delta_i u_{ij} = 0,
\end{equation*}
whence it follows that for any system $ c_j $ the conditions
\begin{equation*}
    \Re \sum\limits_{j=1}^{n} c_j u_{ij} > 0, ~~~ i=1,\dots,m,
\end{equation*}
cannot be satisfied simultaneously.

Sufficiency. Introduce the function
\begin{equation*}
    v = v(\delta_1, \dots, \delta_m) = \sum\limits_{j=1}^{n} \left| \sum\limits_{i=1}^{m} \delta_i u_{ij} \right|, ~~~ \delta_i \ge 0, ~~~ \sum\limits_{i=1}^{m} \delta_i = 1.
\end{equation*}
Since the function $ v $ is continuous on a compact set, it attains the minimum value $ v_0 $ for $ \delta_i = \delta_i^0 $.
Let us show that the condition of the lemma is equivalent to the following: if the inequalities $ \Re \sum \limits_{j = 1}^{n} c_j u_{ij}> 0 $ cannot be satisfied simultaneously, then $ v_0 = 0 $.

Let us prove it by contradiction. Let the inequalities be satisfied simultaneously for any $ c_j $, but $ v_0> 0 $. Take 
\begin{equation*}
    c_j = \sum\limits_{i=1}^{m} \delta_i^0 \overline{u}_{ij}
\end{equation*}
and assume that for such $ c_j $, without loss of generality, the inequality does not hold for $ i = m $
\begin{equation*}
    \Re \sum\limits_{j=1}^{n} \left( \sum\limits_{i=1}^{m} \delta_i^0 \overline{u}_{ij} \right) u_{mj} \le 0.
\end{equation*}
Denote
\begin{equation*}
    v_* = \sum\limits_{j=1}^{m} |u_{mj}|^2, ~~~ \lambda = \dfrac{v_*}{v_* + v_0} < 1.
\end{equation*}
Let us choose $\delta_i$ as
\begin{equation*}
    \delta_i = \begin{cases}
    \lambda \delta_i^0, & i=1,2,\dots,m-1, \\
    (1 - \lambda) + \lambda \delta_m^0, & i = m, \\
    \end{cases}
\end{equation*}
and show that $v(\delta_1, \dots, \delta_m) < v_0$. Indeed,
\begin{eqnarray*}
v & = & \sum\limits_{j=1}^{n} \left| \sum\limits_{i=1}^{m} \delta_i u_{ij} \right| = 
    \sum\limits_{j=1}^{n} \left| (1 - \lambda)u_{mj} + \lambda \sum\limits_{i=1}^{m} \delta_i^0 u_{ij} \right| \\
& = & 
    (1 - \lambda)^2 \sum\limits_{j}^{n} |u_{mj}|^2 + \lambda^2 \sum\limits_{j=1}^{n} \left| \sum\limits_{i=1}^{m} \delta_i^0 u_{ij} \right| + 
    2 \lambda (1 - \lambda) \Re \left( \sum\limits_{j=1}^{n} \sum\limits_{i=1}^{m} \delta_i^0 \overline{u}_{ij} u_{mj} \right).
\end{eqnarray*}
As noted above, 
\begin{equation*}
    \Re \left( \sum\limits_{j=1}^{n} \sum\limits_{i=1}^{m} \delta_i^0 \overline{u}_{ij} u_{mj} \right) \le 0,
\end{equation*}
but $\lambda \ge 0$, $1 - \lambda > 0$, whence
\begin{equation*}
    v \le (1 - \lambda)^2 v_* + \lambda^2 v_0 = \dfrac{v_0^2}{(v_* + v_0)^2} v_* + \dfrac{v_*^2}{(v_* + v_0)^2} v_0 =
    v_0 v_* \dfrac{v_* + v_0}{(v_* + v_0)^2} = \dfrac{v_*}{v_* + v_0} v_0 = \lambda v_0 < v_0.
\end{equation*}
We come to a contradiction with the optimality of $ v_0 $, hence $ v_0 = 0 $. 
\end{proof}
\end{lemma}
{\bf Note 1.}
If all $ u_{ij} \in \mathbb{R} $, then it is sufficient to choose $ c_j $ real.

Using this lemma and Kolmogorov's criterion, we prove another optimality criterion.
\begin{theorem}
(Remez). Let a system of vectors with a matrix $ V \in \mathbb{C}^{n \times r} $ and a vector $ a \in \mathbb{C}^n $ be given, so that vector $a$ should be approximated by a linear combination of columns of $ V $. Let the vector $ V \widehat{u} $ attain the maximum absolute values of the residual in the positions $ E = \{i_1, \dots, i_t \} $. Then $ V \widehat{u} $ is the vector of the best approximation for $ a $ if and only if there are $ \delta_k \ge 0 $, $ k = 1, \dots, t $, not all of which are equal zero, such that 
\begin{equation*}
    \sum\limits_{k \in E} \delta_k \overline{\left( a_{k} - \widehat{u}^T v^{k} \right)} v_j^{k} = 0, ~~~ j = 1, \dots, r.
\end{equation*}
\begin{proof}
Sufficiency. Let the conditions
\begin{equation*}
    \sum\limits_{k \in E} \delta_k \overline{\left( a_{k} - \widehat{u}^T v^{k} \right)} v_j^{k} = 0, ~~~ j = 1, \dots, r.
\end{equation*}
hold. Denote $u_{ij} = \overline{\left( a_{k} - \widehat{u}^T v^{k} \right)} v_j^{k}$ as in the previous lemma. Then, according to the lemma, the conditions
\begin{equation*}
    \Re \sum\limits_{j \in E} c_j v_j^{k} \overline{\left( a_{k} - \widehat{u}^T v^{k} \right)} > 0
\end{equation*}
are not fulfilled simultaneously for all $ c_j $. Taking into account that
$
    \sum\limits_{j \in E} c_j v_j^{k}
$
defines an arbitrary vector of the form $Vc$, we get that the Kolmogorov criterion is fulfilled and $V\widehat{u}$ is the optimal approximation.

Necessity. If $V\widehat{u}$ is optimal, then Kolmogorov's criterion is fulfilled
\begin{equation*}
    \min\limits_{k \in E} \Re \sum\limits_{j \in E} c_j v_j^{k} \overline{\left( a_{k} - \widehat{u}^T v^{k} \right)} \le 0,
\end{equation*}
hence,
\begin{equation*}
    \Re \sum\limits_{j \in E} c_j v_j^{k} \overline{\left( a_{k} - \widehat{u}^T v^{k} \right)} > 0
\end{equation*}
are not fulfilled simultaneously, and then by the lemma there are $\delta_k\ge 0$ satisfying the Remez conditions.
\end{proof}
\end{theorem}
{\bf Note 2.}
There is always an optimal solution in which the maximum values are attained at $t$ points, where $1\le t\le 2r + 1$ in the complex case and $1\le t\le r + 1$ in the real case, and the condition of the Remez criterion is satisfied with $\delta_k > 0$.

{\bf Note 3.}
Note that the Remez condition can be rewritten as
\begin{equation*}
    \sum\limits_{k \in E} \delta_k \sign \left\{\overline{\left( a_{k} - \widehat{u}^T v^{k} \right)} \right\} v_j^{k} = 0, ~~~ j = 1, \dots, r.
\end{equation*}
Suppose we somehow found the set $E$. Then using the Remez criterion it is easy to find a solution. Let us solve the system
\begin{equation*}
    \sum\limits_{k \in E} v_j^{k} s_k = 0, ~~~ j = 1, \dots, r,
\end{equation*}
This is a system with $r + 1$ variables and $r$ equations, which has a non-trivial solution. Let us find this solution. Then
\begin{equation}
\label{equidistant_signs}
    \delta_k = |s_k|, ~~~ \sign \left\{\overline{\left( a_{k} - \widehat{u}^T v^{k} \right)} \right\} = \sign s_k.
\end{equation}
Since Chebyshev system always has a characteristic set consisting exactly of $r+1$ elements, for Chebyshev system with $r+1$ equations and $r$ unknowns, the equation (\ref{equidistant_signs}) gives the signs of the quantities $\overline{\left(a_{k} - \widehat{u}^T v^{k}\right)}$ for the best equidistant point. In addition, the absolute values of $\left|a_{k} - \widehat{u}^T v^{k} \right|$ are equal, which allows in the real case to write out a system with $r$ equations and $r$ unknowns, the solution of which is the best equidistant point. Thus, the best equidistant point of the Chebyshev system of size $(r+1)\times r$ can be found for $O(r^3)$ operations by solving a system of linear equations twice.

%\begin{note}
%Если применить критерий Ремеза к системе из $r + 1$ уравнений и $r$ неизвестных, то найденные знаки будут соответствовать наилучшей равноудаленной точке системы. Отсюда следует, что условие Ремеза для выбранной системы строк просто дает знаки невязок, гарантирующие, что равноудаленная точка с такими знаками будет оптимальна. Следовательно, критерий Ремеза можно еще сформулировать следующим образом. Если есть система из $r+1$ позиции такая, что для оптимального решение на этих строках невязка будет принимать максимальные значения в позициях этих строк, то это оптимальное решение.
%\end{note}

% \section{О задаче поиска характеристических множеств в вещественном случае}
% \label{real_case_theory}
\begin{center}
{\large 4. ON THE PROBLEM OF FINDING CHARACTERISTIC SETS IN THE REAL CASE }
\end{center}
Let us move on to the methods for solving the problem (\ref{vector_problem_statement}).

% \subsection{Комбинаторная формула решения}
\begin{center}
	{\it 4.1. Combinatorial solution formula}
\end{center}

Suppose we need to solve the problem of approximating the vector $a\in\mathbb{R}^{n}$ by the system of vectors $V\in\mathbb{R}^{n\times r}$. Let $\widetilde{V} \in \mathbb{R}^{(r + 1) \times r}$ and $\widetilde{a} \in \mathbb{R}^{r + 1}$. Then we denote by $\begin{bmatrix} 
    \widetilde{V} & \widetilde{a}
\end{bmatrix} \in \mathbb{R}^{(r + 1) \times (r + 1)}$ the matrix whose first $r$ columns are columns of the matrix $\widetilde{V}$, and the last column is $\widetilde{a}$. Denote by $V(i_1, \dots, i_k)$ the submatrix of the matrix $V$ containing the rows $i_1, \dots, i_k$. Similarly, we denote by $a(i_1, \dots, i_k)$ the subvector of the vector $a$ containing elements $i_1, \dots, i_k$. In addition, we denote by $\widetilde{V}_{\setminus k}$ the submatrix of the matrix $\widetilde{V}$, in which the $k$-th row is deleted.

\begin{theorem}
\label{solution_formula}
Suppose we need to solve the problem of approximating the vector $a\in\mathbb{R}^{n}$ by the system of vectors $V\in\mathbb{R}^{n\times r}$. Let
\begin{equation*}
    \mu = \inf\limits_{u \in \mathbb{R}^{r}} \left\| a - Vu \right\|_{\infty}.
\end{equation*}
Then
\begin{equation*}
    \mu = \max\limits_{i_1, i_2, \dots, i_{r + 1}} \dfrac{\left|  \det \begin{bmatrix}
    V(i_1, i_2, \dots, i_{r+1}) & a(i_1, i_2, \dots, i_{r+1})
\end{bmatrix}  \right|}{\sum\limits_{k=1}^{r+1} \left| \det \left( V(i_1, i_2, \dots, i_{r+1})_{\setminus k} \right) \right|}.
\end{equation*}
\begin{proof}
By the Theorem~\ref{schnirelman_matrix} we have that
\begin{equation*}
    \mu(J) = \mu_{r+1}(J),
\end{equation*}
and according to the definition
\begin{equation*}
    \mu_{r+1}(J) = \max\limits_{J_{r+1} \in M_{r+1}} \mu(J_{r+1}).
\end{equation*}
It remains to apply the Theorem~\ref{best_equidistant_solution} to obtain the explicit form of $\mu(J_{r+1})$.

\end{proof}
\end{theorem}

% \subsection{Аналог теоремы о чебышевском альтернансе}
\begin{center}
	{\it 4.2. An analogue of the Chebyshev alternance theorem}
\end{center}

In the statement~\ref{optimal_solution_has_equal_elements}, it was shown that for the optimal solution in the residual vector, the maximum absoulte values are achieved in at least $r+1$ positions. The result on the Chebyshev alternance for continuous functions is widely known. For continuous functions it states, in addition to the fact that the maximum abosulte values are reached in the residual vector, that there is an alternation of signs. A similar result can be proved in the matrix case.

\begin{lemma}
\label{residual_alternating_signs}
Suppose we need to solve the problem of approximating the vector $a\in\mathbb{R}^{n}$ by the system of vectors $V\in\mathbb{R}^{n\times r}$ and that the vector $z^*$ is the vector of the best approximation (the best equidistant point of the system). Denote by $w = a - Vz^*$ the residual vector. Then the signs of the quantities
\begin{equation}
\label{sign_alternance}
    w_1 D_1, w_2 D_2, \dots, w_{n+1} D_{n+1}
\end{equation}
alternate.
\begin{proof}
In the notation of Theorems~\ref{full_equidistant_solution} and~\ref{best_equidistant_solution} we have that
\begin{equation*}
    z^* = \sum\limits_{j=1}^{n+1} \dfrac{|D_j|}{\sum\limits_{\nu=1}^{n+1} |D_\nu|} z^j.
\end{equation*}
Moreover, from the Theorem~\ref{full_equidistant_solution} we have that
\begin{equation*}
    (v^j, z^j) - a_j = \dfrac{(-1)^{j+1}}{D_j} \sum\limits_{\nu=1}^{n+1} (-1)^{\nu} a_{\nu} D_{\nu} = \dfrac{(-1)^{j+1}}{D_j} X,
\end{equation*}
where
\begin{equation*}
    X = \sum\limits_{\nu=1}^{n+1} (-1)^{\nu} a_{\nu} D_{\nu}
\end{equation*}
does not depend on $j$.
According to the definition $z^j$,
\begin{equation*}
    Vz_j = \begin{bmatrix}
        a_1 &
        a_2 &
        \dots &
        a_{j-1} &
        \widetilde{a}_j &
        a_{j+1} &
        \dots &
        a_{n+1}
    \end{bmatrix}^T,
\end{equation*}
where $\widetilde{a}_j = (v^j, z^j)$. Then
\begin{equation*}
    Vz^* = \sum\limits_{j=1}^{n+1} \dfrac{|D_j|}{\sum\limits_{\nu=1}^{n+1} |D_\nu|} \begin{bmatrix}
        a_1 \\
        a_2 \\
        \vdots \\
        a_{j-1} \\
        \widetilde{a}_j \\
        a_{j+1} \\
        \vdots \\
        a_{n+1}
    \end{bmatrix} = \begin{bmatrix}
        a_1 \\
        a_2 \\
        \vdots \\
        a_{j-1} \\
        a_j \\
        a_{j+1} \\
        \vdots \\
        a_{n+1}
    \end{bmatrix} - \dfrac{1}{\sum\limits_{\nu=1}^{n+1} |D_\nu|} \left(
    \begin{bmatrix}
        |D_1| \widetilde{a}_1 \\
        |D_2| \widetilde{a}_2 \\
        \vdots \\
        |D_{j-1}| \widetilde{a}_{j-1} \\
        |D_j| \widetilde{a}_j \\
        |D_{j+1}| \widetilde{a}_{j+1} \\
        \vdots \\
        |D_{n+1}| \widetilde{a}_{n+1}
    \end{bmatrix} - \begin{bmatrix}
        |D_1| a_1 \\
        |D_2| a_2 \\
        \vdots \\
        |D_{j-1}| a_{j-1} \\
        |D_j| a_j \\
        |D_{j+1}| a_{j+1} \\
        \vdots \\
        |D_{n+1}| a_{n+1}
    \end{bmatrix}
    \right).
\end{equation*}
From here
\begin{equation*}
    a - Vz^* = \dfrac{1}{\sum\limits_{\nu=1}^{n+1} |D_\nu|} \cdot \begin{bmatrix}
        |D_1| (a_1 - \widetilde{a}_1) \\
        |D_2| (a_2 - \widetilde{a}_2) \\
        \vdots \\
        |D_{n+1}| (a_{n+1} - \widetilde{a}_{n+1}) \\
    \end{bmatrix} = \dfrac{1}{\sum\limits_{\nu=1}^{n+1} |D_\nu|} \cdot \begin{bmatrix}
        |D_1| \dfrac{(-1)^1}{D_1} X \\
        |D_2| \dfrac{(-1)^2}{D_2} X \\
        \vdots \\
        |D_{n+1}| \dfrac{(-1)^{n+1}}{D_{n+1}} X \\
    \end{bmatrix},
\end{equation*}
since
\begin{equation*}
    a_j - \widetilde{a}_j = a_j - (v^j, z^j) = \dfrac{(-1)^j}{D_j} X.
\end{equation*}
Denoting
\begin{equation*}
    C = \dfrac{X}{\sum\limits_{\nu=1}^{n+1} |D_\nu|},
\end{equation*}
we get that
\begin{equation*}
    w = a - Vz^* = C \begin{bmatrix}
        (-1)^{1} \sign D_{1} \\
        (-1)^{2} \sign D_{2} \\
        \vdots \\
        (-1)^{n+1} \sign D_{n+1} \\
    \end{bmatrix}.
\end{equation*}
And then
\begin{equation*}
    w_j D_j = C (-1)^j |D_j|,
\end{equation*}
and the sequence
\begin{equation*}
    w_1 D_1, w_2 D_2, \dots, w_{n+1} D_{n+1}
\end{equation*}
has alternating signs.
\end{proof}
\end{lemma}

% \subsection{Обобщенный алгоритм Ремеза для матриц}
\begin{center}
	{\it 4.3. Generalized Remez algorithm for matrices}
\end{center}

Note that Lemma~\ref{residual_alternating_signs} shows that the matrix problem of the best approximation is more general than a similar problem for continuous functions. Indeed, for continuous functions, as well as for matrices, the result is known that there exists a characteristic set of $r+1$ elements when approximating by a system of $r$ functions (i.e. by polynomials of degree $r-1$). With a known characteristic set, the problem is reduced to solving a matrix problem with a Vandermond matrix. Note that the determinant of the Vandermonde matrix can be calculated by the formula
\begin{equation*}
    W(x_1, x_2, \dots, x_r) = \prod\limits_{j<i} (x_i - x_j).
\end{equation*}
Hence it is easy to see that the sign of the determinant of the Vandermonde matrix depends only on the order in which the points are taken. If the points are taken in ascending order each time, then all determinants have the same sign and in the formula (\ref{sign_alternance}) only the signs of the residual elements remain. This reasoning allows us to generalize the Remez algorithm for constructing the best Chebyshev approximation to the matrix case.

Note that the theorem~\ref{solution_formula} already allows us to solve the problem (\ref{vector_problem_statement}) in a finite number of operations. To do this, it is enough to iterate over all possible variants of characteristic sets (all possible sets of $r+1$ rows) and solve for each of them the problem of finding the best equidistant point. However, it is possible to find the characteristic set much faster.

Let us present an algorithm for solving the problem of the best Chebyshev approximation in the real case. For given matrix $ V \in \mathbb{R}^{n \times r} $ and vector $ a \in \mathbb{R}^{n} $:
\begin{enumerate}
	\item Choose an arbitrary set of $r+1$ indices of the rows of the matrix $V$. Denote this set by $I_1$ and take $t=1$.
	\item Solve the problem of the best uniform approximation for the matrix $V(I_t)$ and the vector $a(I_t)$. This problem can be solved in $O(r^3)$ operations for the Chebyshev system of vectors. Denote the solution by $u_t$.
	\item Calculate the residual $ w_t = Vu_t - a $ and find the element of $w_t$ with the maximum absolute value. This requires $ O (nr) $ operations. Denote the position of this element by $ j_t $. If $ j_t \in I_t $, then according to the note to the Remez criterion, the set $ I_t $ is characteristic and $ u_t $ is the solution to the problem of the best approximation.
	\item If $ j_t \notin I_t $, then we try to replace each of the elements of the set $ I_t $ by $ j_t $. Let $ I_t = \{i_1^t, i_2^t, \dots, i_{r + 1}^t \} $. We denote $ I_t^k = I_t \setminus \{i_k^t \} \cup j_t $. Solve the problem with the matrix $ V(I_t^k) $ and the vector $ a(I_t^k) $ and find the maximum of the absolute values in the residual on the set $ I_t^k $, $ w_t^k = V(I_t^k) u_t^k - a(I_t^k) $. Let $ l = \argmax \limits_{k} \| w_t^k \|_\infty $. This step requires $ O (r ^ 4) $ operations 
	\item $I_{t+1} = I_t^l$, $t = t + 1$ and go to step 2.
\end{enumerate}

\begin{theorem}
\label{generalized_remez}
Let a system of vectors $ V \in \mathbb{R}^{n \times r} $ be Chebyshev and let $ a \in \mathbb{R}^{n} $. Then the generalized Remez algorithm finds the solution to the problem of the best approximation in a finite number of operations.
\begin{proof}
Let $ I_t $ be the current set of indices and $ w_t = Vu_t - a $.  Denote
\begin{equation*}
E_t = \| w_t(I_t) \|_\infty.
\end{equation*}
Assume that in the vector $ w_t $ the element with the maximum absolute value is attained at the position $ j_t $. Then consider the problem of the best approximation for the submatrix taken on the set of rows with numbers $ I_t \cup \{j_t \} $. For this problem, there is a characteristic set of $ r + 1 $ elements. Note that it cannot entirely consist of elements of the set $ I_t $, since for an optimal solution on this set, a strictly larger value of the residual element is attained in the position $ j_t $. This means that the characteristic set contains $ j_t $ and $ r $ elements from the set $ I_t $, that is, it is obtained by replacing one of the elements in $ I_t $ with $ j_t $. Let us denote this set by $ \widehat{I}_t $. Let us show that in this case the error of the optimal approximation on the new set $ \widehat{I}_t $ is strictly greater than on the set $ I_t $.
Indeed, assume that the error of the optimal approximation on the set $ I_t $ is $ \delta $, and that the corresponding solution gives the error $ \varepsilon $ on the element in the position $ j_t $. Note that $ \varepsilon> \delta $. Similarly, let the optimal approximation on the set $ \widehat{I} _t $ on the set itself get the error $ \delta_1 $, and on the element removed from $ I_t $, the error $ \varepsilon_1 $. Note that since $ \widehat {I} _t $ is a characteristic set, then $ \varepsilon \le \delta_1 $. Then
\begin{equation*}
    \delta < \varepsilon \le \delta_1,
\end{equation*}
whence it follows that the error of the optimal approximation on the new set $ \widehat{I}_t $ is strictly greater than the error of the optimal approximation on the set $ I_t $. Hence it follows that 
\begin{equation*}
    E_{t+1} > E_{t}.
\end{equation*}
But since there are finitely many subsets of size $ r + 1 $, the sequence $ \{E_t \} $ cannot be infinite and reaches its maximum value on some set, which, according to the reasoning of Theorem~\ref{solution_formula}, indicates that the found set is characteristic and optimal solution is constructed.

\end{proof}
\end{theorem}

% \section{Задача о чебышевском приближении матриц}
% \label{matrix_case}
\begin{center}
{\large 5. ON THE PROBLEM OF THE CHEBYSHEV APPROXIMATION OF MATRICES }
\end{center}
Having built a method for solving the problem (\ref{one_matrix_problem_statement}), we can go to the problem
\begin{equation}
\label{global_problem_statement}
\mu = \inf\limits_{U \in \mathbb{R}^{m \times r}, V \in \mathbb{R}^{n \times r}} \left\| A - UV^T \right\|_{C}.
\end{equation}

% \subsection{Необходимое условие оптимальности}
\begin{center}
	{\it 5.1. Necessary condition for optimality}
\end{center}

In \cite{daugavet1971}, the necessary condition for the optimality of the solution of the problem (\ref{global_problem_statement}) was proved. In light of the results obtained above, we can easily obtain this condition and also generalize it to the case of arbitrary rank. Let the pair $ (\widehat{U}, \widehat{V}) $ be a solution to the problem (\ref{global_problem_statement}). Then the matrix $ \widehat {U} $ is a solution to the problem 
\begin{equation*}
\mu = \inf\limits_{U \in \mathbb{R}^{m \times r}} \left\| A - U\widehat{V}^T \right\|_{C}.
\end{equation*}
Assume that that the element of the matrix $ A - \widehat{U} \widehat{V}^T $ with maximum absolute value is at position $ (i, j) $. Consider the problem for the $i$-th row of the matrix 
\begin{equation}
\label{vector_subproblem}
\mu = \inf\limits_{u \in \mathbb{R}^{m \times r}} \left\| a^i - u \widehat{V}^T \right\|_{C}.
\end{equation}
It is clear that $ \widehat{u}^i $ is the optimal solution to the problem (\ref{vector_subproblem}), otherwise it would be possible to replace the $i$-th row of the matrix $ \widehat{U} $ with the optimal solution and obtain better result in problem (\ref{global_problem_statement}).
By virtue of the optimality of the solution, we have that in the vector $ a^i - \widehat{u}^i \widehat{V}^T $ the maximum absolute value is attained at $ r + 1 $ position and the signs of the residual and determinants of the matrix $ \widehat {V} $ in these positions alternate (see Lemma~\ref{residual_alternating_signs}). To obtain the necessary condition from \cite{daugavet1971}, it suffices to note that for $ r = 1 $ the determinants are the elements of the vector $ \widehat{V} $, and the sign of $ \widehat{u}^i $ obviously does not change within one column from which Statement~\ref{simple_necessary_condition} follows. In the case of arbitrary rank in each column and each row in which there is an element with maximum absolute value, the maximum absolute value is reached at $ r + 1 $ positions and the signs of the residual and determinants alternate according to Lemma~\ref{residual_alternating_signs}.

% \subsection{Метод решения}
% \label{solution_method}
\begin{center}
	{\it 5.2. Solving method}
\end{center}

Let us build an iterative process for solving the problem (\ref{global_problem_statement}). Let a matrix $ A \in \mathbb{R}^{m \times n} $ and a Chebyshev matrix $ U_0 $ be given. Let us find the best approximation $ A = U_0 V ^ T $ and denote the result by $ V_1 $. Suppose that the system of vectors $ V_1 $ is Chebyshev. Then we find the best approximation $ A = U V_1 ^ T $ and denote the result by $ U_1 $, again assuming that the system $ U_1 $ is Chebyshev. Let us find the best approximation $ A = U_1 V ^ T $ and denote the result by $ V_2 $. Continuing according to the described scheme, we note that the quantity $ \rho_k = \| A - U_k V_k ^ T \|_{C} $ does not increase and is bounded below and, therefore, converges.

% \section{Численные эксперименты}
% \label{num_exp}
\begin{center}
{\large 6. NUMERICAL EXPERIMENTS}
\end{center}
In this section, we present a number of numerical experiments, where we apply the method described in Section~5.2 to construct low-rank Chebyshev approximations for matrices whose singular values do not decrease. For the experiments, the algorithm from Section~5.2 was implemented in C++. In the experiment, random matrices were generated with singular values uniformly distributed on $ [1, 2] $. For this, two random matrices are generated from the standard normal distribution, a QR decomposition is constructed for them, and the factors $ Q $ are chosen as the left and right singular vectors. 
The singular value matrix $ \Sigma $ is generated as a diagonal matrix such that its diagonal elements are uniformly distributed on $ [1, 2] $. After that, the matrix $ U \Sigma V^T $ is built. The sizes of the matrices vary from 10 to 1400 with a step of 10, and the approximation rank is chosen as $ r = \sqrt{n} $, where $ n $ is the size of the matrix.
For each size, 10 random matrices are generated, and for each of them, the alternance method of 20 random points is launched. Thus, for each size,
there are 200 runs of the alternance method. For each matrix, the mean value of the accuracy $ \mu_n^i $ and the variance $ \sigma_n^i $ of the Chebyshev approximation over 20 starting points are calculated. Further, for each size, these values are averaged over 10 matrices
\begin{equation*}
    \mu_n = \dfrac{1}{10} \sum\limits_{i=1}^{10} \mu_n^i, ~~~ \sigma_n = \dfrac{1}{10} \sum\limits_{i=1}^{10} \sigma_n^i.
\end{equation*}

\begin{figure}%
    \centering
    \subfloat[\centering Approximation error]{{\label{figure_precision}\includegraphics[width=93mm]{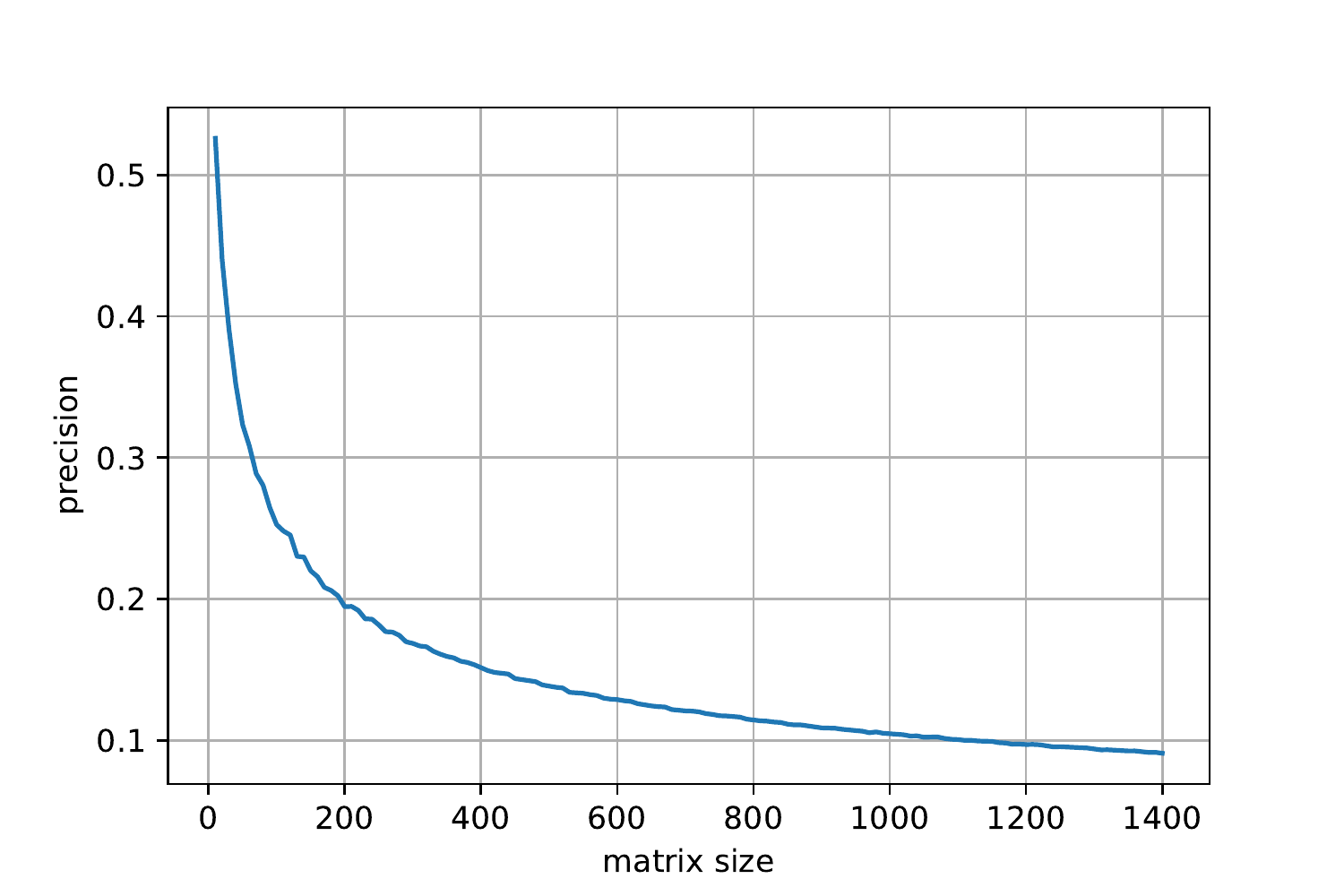} }}%
    % \qquad
    \subfloat[\centering Variance]{{\label{figure_variation}\includegraphics[width=93mm]{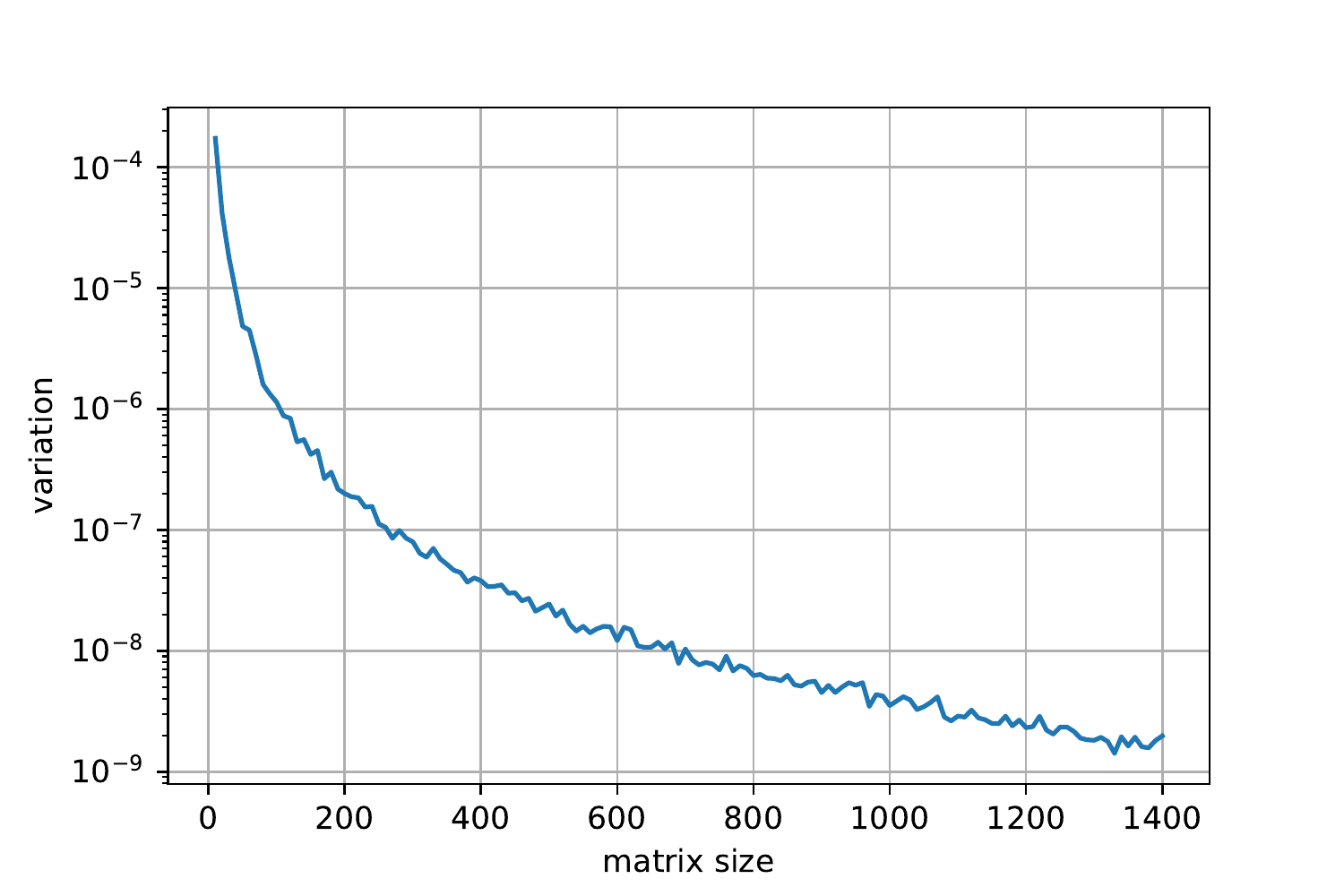} }}%
    
    Fig. 1. Averaged approximation error and variance over 20 random initial conditions and 10 random matrices for different matrix sizes.%
\end{figure}

\begin{figure}[!ht]
\centering     %%% not \center
\includegraphics[width=120mm]{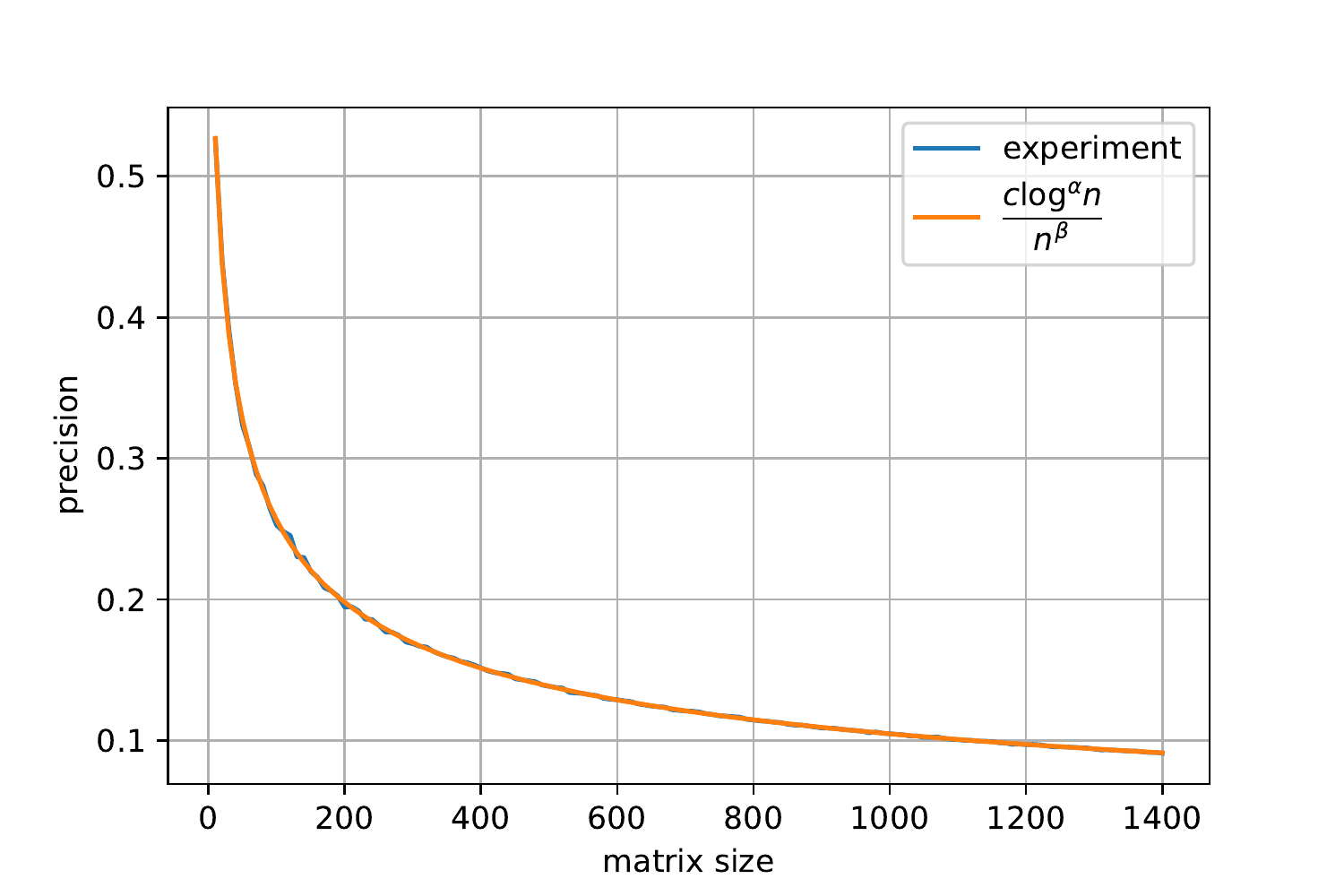}

Fig. 2. Averaged approximation error over 20 random initial conditions and 10 random matrices of different sizes, and its fitted approximation.
\label{figure_precision_with_function}
\end{figure}

\begin{figure}[!ht]
\centering     %%% not \center
\includegraphics[width=120mm]{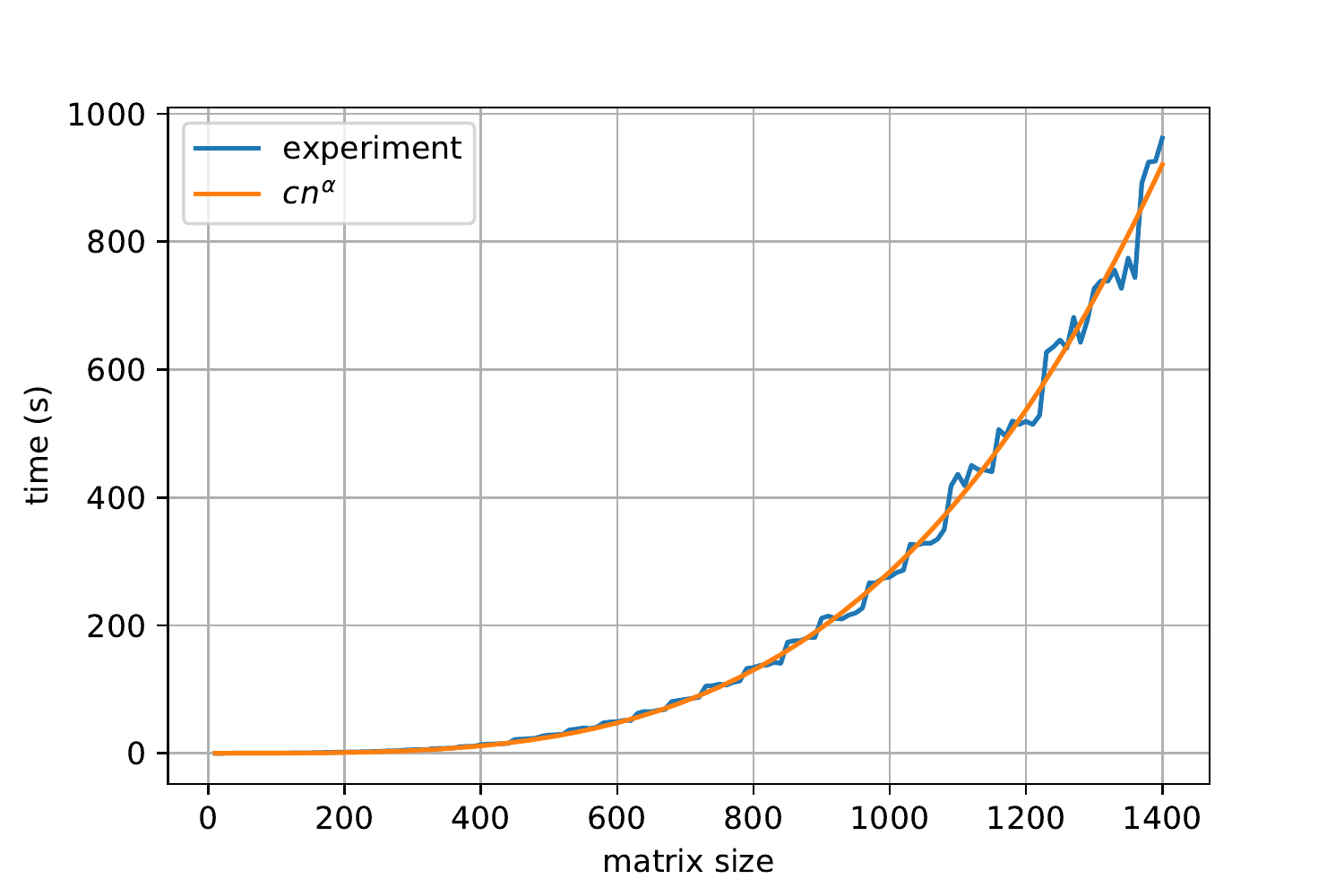}

Fig. 3. Averaged execution time over 200 runs for matrices of different sizes, and its fitted approximation.
\label{figure_time_with_function}
\end{figure}

In Fig.~1a, we plot $\mu_n$ against the size of the matrix, and in Fig.~1b we present the corresponding plot for $\sigma_n$; the latter is shown in logarithmic scale. An interesting observation is that variance decays as the size of the problem grows. For instance, for a matrix of sizes $1400 \times 1400$, if we approximate it with rank $37$ with 20 different random initial conditions, we get the following approximation errors:
\begin{equation*}
    \begin{bmatrix}
       0.09111796 & 0.09098914 & 0.09103979 & 0.09101653 & 0.09097955 \\
       0.09112523 & 0.09102086 & 0.09097652 & 0.09099676 & 0.0909908  \\
       0.09106326 & 0.0911168  & 0.09108753 & 0.09101277 & 0.09098213 \\
       0.09103401 & 0.09106984 & 0.09097417 & 0.09094869 & 0.09092307
    \end{bmatrix}.
\end{equation*}
We also estimated the asymptotic dependence of the approximation error on the matrix sizes. To this end, we fitted a parametric curve $\dfrac{c \log^{\alpha} n}{n^\beta}$. The optimal parameter values that we found are
\begin{equation*}
    \begin{cases}
    c = 0.995139, \\
    \alpha = 0.604346, \\
    \beta = 0.495001. \\
    \end{cases}
\end{equation*}
The curve $\dfrac{c \log^{\alpha} n}{n^\beta}$ with the optimal parameters is presented in Fig.~2. It is worth noting that this estimate corresponds to
\begin{equation*}
    \varepsilon \approx \dfrac{\log^{0.6}{n}}{n^{0.5}},
\end{equation*}
while the known proved bound (see Thm.~~\ref{townsend_theorem} \cite{udell2019big}) for $r = \sqrt{n}$ gives
\begin{equation*}
    \varepsilon \le \dfrac{6\sqrt{2} \log^{0.5}{(2n + 1)}}{n^{0.25}},
\end{equation*}
which shows that the latter is not optimal.

In addition, we measured the execution time of our algorithm. For every problem size, we averaged the time over 200 runs. In Fig.~3, we plot the execution time against the matrix size. We also approximated the experimental curve by fitting $c n^\alpha$. The optimal parameters are
\begin{equation*}
    \begin{cases}
    c = 9.13302\mathrm{e}{-09} \\
    \alpha = 3.49745 \\
    \end{cases}
\end{equation*}
This means that for $r = \sqrt{n}$ the computational complexity is $O(n^{3.5})$ in practice.

% \section{Заключение}
% \label{conclusion}
\begin{center}
{\large 7. CONCLUSION}
\end{center}
In this work, a method is proposed for solving the problem of the best low-rank approximation of matrices in the Chebyshev norm in the case when one of the factors of skeleton decomposition is known. Using the proposed scheme and the alternance method, an algorithm was constructed for computing low-rank approximations in the Chebyshev norm for an arbitrary rank. The described method essentially generalizes all known methods for solving the problem of Chebyshev approximations of matrices and improves the known theoretical estimates of the approximation described in \cite{udell2019big}. Numerical experiments show that the method is capable of approximating matrices with good accuracy even in the absence of decay of singular values and has an acceptable asymptotic complexity.

\clearpage
\begin{center}

\end{center}

\end{document}